\documentclass[12pt]{article}

\usepackage{amscd,amsmath, amssymb, fancyhdr, mathbbol}
\numberwithin{equation}{section}

\newcommand{\version}{version 1.2,\ \ January 21, 2018}

\def\eqref#1{(\ref{#1})}

\newcommand{\arrow}{{\:\longrightarrow\:}}
\newcommand{\Z}{{\mathbb{Z}}}
\def\C{{\Bbb C}}

\newcommand{\calo}{{\mathcal{O}}}
\newcommand{\restrict}[1]{{\left|_{{\phantom{|}\!\!}_{#1}}\right.}}

\newcommand{\Ext}{\operatorname{Ext}}
\newcommand{\Hom}{\operatorname{Hom}}

\newcounter{Mycounter}[section]
\newcounter{lemma}[section]
\newcounter{claim}[section]
\newcounter{corollary}[section]
\renewcommand{\thecorollary}{{Corollary \thesection.\arabic{corollary}}}
\newcommand{\corollary}{%
    \setcounter{corollary}{\value{Mycounter}}
    \refstepcounter{corollary}
    \stepcounter{Mycounter}
    {\noindent \bf \thecorollary:\ }}

\newcounter{theorem}[section]
\renewcommand{\thetheorem}{{Theorem \thesection.\arabic{theorem}}}
\newcommand{\theorem}{%
    \setcounter{theorem}{\value{Mycounter}}
    \refstepcounter{theorem}
    \stepcounter{Mycounter}
    {\noindent \bf \thetheorem:\ }}

\newcounter{conjecture}[section]
\newcounter{proposition}[section]
\renewcommand{\theproposition}
      {{Proposition \thesection.\arabic{proposition}}}
\newcommand{\proposition}{%
    \setcounter{proposition}{\value{Mycounter}}
    \refstepcounter{proposition}
    \stepcounter{Mycounter}
    {\noindent \bf \theproposition:\ }}

\newcounter{definition}[section]
\renewcommand{\thedefinition}
      {{Definition~\thesection.\arabic{definition}}}
\newcommand{\definition}{%
    \setcounter{definition}{\value{Mycounter}}
    \refstepcounter{definition}
    \stepcounter{Mycounter}
    {\noindent \bf \thedefinition:\ }}

\newcounter{example}[section]
\newcounter{remark}[section]
\newcommand{\proof}{\noindent{\bf Proof:\ }}

\def\blacksquare{\hbox{\vrule width 5pt height 5pt depth 0pt}}
\def\endproof{\blacksquare}

\makeatletter
\@addtoreset{equation}{section} \@addtoreset{footnote}{section}
\makeatother

\begin{document}
\begin{center}
{\LARGE\bf
Multiplicity of singularities is not a bi-Lipschitz invariant\\[4mm]
}

Lev Birbrair\footnote{Partially supported 
by  CNPq grant 302655/2014-0.}, Alexandre Fernandes\footnote{Partially supported 
by  CNPq grant 302764/2014-7.}, J. Edson Sampaio,
Misha Verbitsky\footnote{Partially supported 
by  the  Russian Academic Excellence Project '5-100'.}
\end{center}

{\small \hspace{0.10\linewidth}
\begin{minipage}[t]{0.85\linewidth}
{\bf Abstract.} 
It was conjectured that multiplicity of a singularity is
bi-Lipschitz invariant. We disprove this conjecture, 
constructing examples of bi-Lipschitz equivalent complex
algebraic singularities with different values of multiplicity.
\end{minipage}
}

\tableofcontents


\section{Introduction}

The famous multiplicity conjecture, stated by Zariski in 1971 (see \cite{Zariski:1971}), is formulated as follows: if two germs of complex analytic hypersurfaces are ambient topolological equivalent, then they have the same multiplicity. It was proved by Zariski \cite{Zariski:1932} for germs of plane analytic curves. The results of Pham-Teisser in \cite{P-T} show that this result can be extended in the following "metric" way: if the two germs of complex analytic curves in n-dimensional space are bi-Lipschitz equivalent with respect to the outer metric, then the germs of the space curves have the same multiplicity. Comte in \cite{Comte:1998} proves that the multiplicity of complex analytic germs (not necessarily codimension 1 sets) is invariant under bi-Lipschitz homeomorphism with Lipschitz constants close enough to 1 (this is a severe assumption). These results motivated the following  question, closely related to the multiplicity conjecture: is the multiplicity of a germ of analytic singularity a bi-Lipschitz invariant? This question was stated as a conjecture in \cite{_BFS}.

 The Lipschitz Regularity Theorem in \cite{BFLS} shows
 that if the multiplicity of a complex analytic germ is
 equal to one, then it is a bi-Lipschitz
 invariant. Namely, if a germs of an analytic set is
 bi-Lipschitz equivalent to a smooth germ, then it is
 smooth itself. Later, Fernandes and Sampaio in
 \cite{FernandesS:2016} give a positive answer to this
 question for surfaces in 3-dimensional space with respect
 to the ambient bi-Lipschitz equivalence. More recently,
 Neumann and Pichon (\cite{NeumannP:2016}) showed that the
 multiplicity is an invariant under bi-Lipschitz
 equivalence, for germs of normal surface singularities.
 Another important result in \cite{FernandesS:2016} is the
 following: in order to prove (or disprove) the
 bi-Lipschitz invariance of the multiplicity, it is enough
 to prove it for the algebraic cones, i.e. for the
 algebraic sets, defined by homogeneous polynomials. In
 \cite{_BFS} the authors show that the conjecture has a
 positive answer for 1 or 2 dimensional complex analytic
 sets.

 The present paper shows that the multiplicity of complex algebraic sets is not a bi-Lipschitz invariant for the sets of dimension bigger or equal to three. Moreover, we show that there exists an infinite family of $3-$dimensional germs, such that all of them are bi-Lipschitz equivalent, but they have different multiplicities. 
The idea of the construction is to consider the complex cones over different embeddings of $\C P^1 \times \C P^1$ to complex projective spaces. Using the topology of Smale-Barden manifolds, we show that all the links of such singularities are diffeomorphic. 
That is why the germs of the corresponding cones are bi-Lipschitz equivalent. From the other hand, the multiplicities of the cones at the origin may be explicitly calculated in terms of the embeddings.


\section{Smale-Barden manifolds}


The classification of 5-manifolds is due to S. Smale
(\cite{_Smale:5_}) and D. Barden (\cite{_Barden:5_}).

\hfill

\definition
A simply connected, compact, oriented 5-manifold is
called {\bf Smale-Barden manifold}. 

\hfill

The Smale-Barden manifolds are uniquely determined by their
second Stiefel-Whitney class and the linking form.

\hfill

\theorem (\cite{_Barden:5_})
Let $X, X'$ be two Smale-Barden manifolds.
Assume that $H^2(X)=H^2(X')$ and this isomorphism
is compatible with the linking form and preserves the
second Stiefel-Whitney class. Then $X$ is diffeomorphic to $X'$.
\endproof

\hfill

\corollary 
There exists only two Smale-Barden manifolds $M$ with $H^2(M)=\Z$,
the product $S^2\times S^3$ and the total space of a 
non-trivial $S^3$-bundle over $S^2$ (see \cite{_Crowley:manifolds_}
for an introduction to Barden theory, where this
manifold is formally introduced).

\hfill

\proof
Indeed, the linking form on $\Z$ vanishes, therefore
the manifold is uniquely determined by the Stiefel-Whitney
class $w_2(M)$. Hence we have only two possibilities:
$w_2(M)=0$ and $w_2(M)\neq 0$. \endproof

\hfill

In early 2000-ies, the classification of 5-manifolds
attracted interest coming from algebraic geometry,
in the context of Sasakian geometry and geometry
of generalized Seifert manifolds (\cite{_Kollar:Seifert_}, 
\cite{_Kollar:circle_}). 
In the present paper we are interested in $S^1$-bundles
over $\C P^1\times \C P^1$.

\hfill

\proposition\label{_S^1_over_CP^1^2_Proposition_}
Let $\pi:\; M\arrow B$ be a simply connected 5-manifold obtained
as a total space of an $S^1$-bundle $L$ over $B=\C P^1\times \C P^1$.
Then $H^2(M)$ is torsion-free, and $M$ is diffeomorphic to $S^2\times S^3$.

\hfill

{\bf Proof. Step 1:}
Universal coefficients formula gives
an exact sequence
\[
  0 \to \Ext_\Z^1(H_1(M; \Z), \Z) \to H^2(M; \Z)\to \Hom_\Z(H_2(M; \Z), \Z)\to 0.
\]
This implies that $ H^2(M; \Z)$ is torsion-free.

\hfill

{\bf Step 2:} 
Consider the following exact sequence of homotopy groups
\[
0\arrow \pi_2(M)\arrow \pi_2(B)\stackrel \phi \arrow
\pi_1(S^1) \arrow \pi_1(M) \arrow 0
\]
Since $\pi_1(M)=0$, the map $\phi$, representing 
the first Chern class of $L$, is surjective.
This exact sequence becomes 
\[ 0 \arrow \pi_2(M)\arrow \Z^2 \arrow \Z \arrow 0
\]
giving $\pi_2(M)=\Z$, and $H^2(M)=\Z$ because
$H^2(M)$ is torsion-free.

\hfill

{\bf Step 3:}  To deduce \ref{_S^1_over_CP^1^2_Proposition_}
from the Smale-Barden classification, it remains to show that
$w_2(M)=0$. However, $w_2(M)=\pi^*(\omega_2(B))$
(\cite[Lemma 36]{_Kollar:circle_}), and the latter
clearly vanishes, because $w_2(S^2)=0$.
\endproof

\section{Multiplicity of homogeneous singularities}

Given a projective variety $X\subset \C P^n$,
the {\bf projective cone} of $X$ is the union
of all 1-dimensional subspaces $l\subset \C^{n+1}$
such that $l$, interpreted as a point in $\C P^n$,
belongs to $X$. The {\bf link} of $C(X)$
is an intersection of $C(X)$ with a unit
sphere $S^{2n+1}\subset \C^{n+1}$.

Let $A$ be complex algebraic set of $\C^{n+1}$ and $x\in A$. The {\bf multiplicity of} $A$ {\bf at} $x$, denoted by $mult(A,x)$, is defined to be the multiplicity of the maximal ideal of the local ring ${\calo}_{A,x}$. Given a projective variety $X\subset\C P^n$, we see that multiplicity of the projective cone $C(X)$ at the origin $0\in\C^{n+1}$ coincides with degree of $X$ (see \cite{Chirka:1989}, subsection 11.3).

Next, we shall be interested in the following
geometric situation. Let $X\subset \C P^n$ be a
variety isomorphic to $\C P^1 \times \C P^1$. Then the Picard group of $X=\C P^1 \times \C P^1$
is isomorphic to $\Z^2$, and a line bundle is
determined by its {\em bidegree}. We shall denote a line
bundle of bidegree $a, b$ by $\calo(a,b)$.

\hfill

\proposition\label{_CP^1^2_link_Proposition_}
Let $X\subset \C P^n$ be a
variety isomorphic to $\C P^1 \times \C P^1$,
and $S\subset C(X)$ be its link.
Assume that and $\calo(1)\restrict X= \calo(a,b)$. Then 
$X$ has degree $2ab$. If, in addition,
$a$ and $b$ are relatively prime, the
link of $C(X)$ is diffeomorphic to $S^2\times S^3$.

\hfill

{\bf Proof. Step 1:}
Since $c_1(\calo(a,b))^2=2ab$, and degree
of a subvariety $X\subset \C P^n$ is its intersection with the
top power of $\calo(1)\restrict X$, one has
$\deg X = 2ab$.

\hfill

{\bf Step 2:} 
Consider the homotopy exact sequence 
\[
0\arrow \pi_2(S)\arrow \pi_2(X)\stackrel \phi \arrow
\pi_1(S^1) \arrow \pi_1(S) \arrow \pi_1(X)\arrow 0
\]
for the circle bundle $\pi:\; S \arrow X$.
Since the map $\phi$ represents the first Chern
class of $\calo(1)\restrict X$, it is obtained
as a quotient of $\Z^2$ by a subgroup generated by
$(a,b)$, and this map is surjective
because $a$ and $b$ are relatively prime.
Then $\pi_1(S)=\pi_1(X)=0$, and 
\ref{_S^1_over_CP^1^2_Proposition_}
implies that $S$ is diffeomorphic to $S^2\times S^3$.
\endproof

\section{Lipschitz invariance of singularities}

Let $X\subset \C^n$ be  a complex variety. The induced metric 
from the Euclidean distance on $\C^n$ gives a distance on $X$; it
is called {\bf the outer metric} on $X$.

\hfill

\definition
Let $X\subset \C^{n}$ and $X'\subset \C^{n'}$ be complex varieties equipped with the outer metrics,
$x\in X, x'\in X'$ marked points.
We say that $(X,x)$ is {\bf bi-Lipschitz equivalent}
to $(X',x')$ if there exist a neighborhoods $U$ of $x$ in $\C^n$ and
$U'$ of $x$ in $\C^n$, and a bi-Lipschitz homeomorphism
of $X\cap U$ to $X'\cap U'$ mapping $x$ to $x'$.
\hfill

\definition
Let $X,X'\subset \C^{n}$ be complex varieties equipped with the outer metrics,
$x\in X, x'\in X'$ marked points.
We say that $(X,x)$ is {\bf ambient bi-Lipschitz equivalent}
to $(X',x')$ if there exists a bi-Lipschitz equivalence
of a neighbourhood $U$ of $x$ in $\C^n$ and a neighbourhood
$U'$ of $x$ in $\C^n$ mapping  $X\cap U$ to $X'\cap U'$  and  $x$ to $x'$.

\hfill

Actually, the two definitions above do not coincide. The
ambient bi-Lipschitz equivalence implies bi-Lipschitz
equivalence, but the examples presented in \cite{BG} show
that the converse does not hold true in general.

As it was already mentioned in Introduction, it was conjectured in \cite{_BFS} that the multiplicity is 
a bi-Lipschitz invariant. We prove that this is false.
Here is the main result of this paper. 

\hfill

\theorem\label{_Main_example_CP^1_times_CP^1_Theorem_}
For each $n\geq 3$, there exists a family $\{ Y_i \}_{i\in\Z}$  of $n$-dimensional complex algebraic varieties
$Y_i \subset \C^{n_i+1}$ such that:
\begin{itemize}
\item [(a)] for each pair $i\neq j$, the germs at the origin of $Y_i\subset \C^{n_i+1}$ and $Y_j\subset \C^{n_j+1}$ are bi-Lipschitz equivalent, but $(Y_i,0)$ and $(Y_j,0)$ have different multiplicity.
\item [(b)] for each pair $i\neq j$, there are $n$-dimensional complex algebraic varieties $Z_{ij}, \widetilde Z_{ij}\subset \C^{n_i+n_j+2}$ such that $(Z_{ij},0)$ and $(\widetilde Z_{ij},0)$ are ambient bi-Lipschitz equivalent, but  $mult(Z_{ij},0)=mult(Y_{i},0)$ and $mult(\widetilde Z_{ij},0)=mult(Y_{j},0)$ and, in particular, they have different multiplicity.
\end{itemize}

\hfill

{\bf Proof.} Let $\{ p_i \}_{i\in\Z}$ be the family of odd prime numbers. 
For each $i\in\Z$, let $L_i$  be a very ample bundle
on $X=\C P^1 \times \C P^1$ of bidegree $(2, p_i)$. Let $X_i$ be projective variety obtained by the embedding of the very ample bundle $L_i$. Consider the 
link of the singularity $S_i:= C(X_i)\cap S^{2m_i+1}$, where $S^{2m_i+1}$ 
is the unit sphere centered in $0\in\C^{m_i+1}$. 
Then, for each pair $i\neq j$ the links $S_i, S_j$ are 
diffeomorphic to $S^2\times S^3$ (\ref{_CP^1^2_link_Proposition_}). 
In particular, $S_i$ to $S_j$ are bi-Lipschitz homeomorphic. Since
a bi-Lipschitz map from $S_i$ to $S_j$ induces a bi-Lipschitz map of their
cones, then the affine cones $(C(X_i),0)$ and $(C(X_j),0)$ are bi-Lipschitz 
equivalent, but $mult(C(X_i),0)=4p_i$ and $mult(C(X_j),0)=4p_j$ (\ref{_CP^1^2_link_Proposition_}).
Thus, if for each $i\in \Z$ we define $Y_i:=C(X_i)\times \C^{n-3}$, then we have that the 
family $\{ Y_i \}_{i\in\Z}$ satisfies the item (a), since 
$mult(Y_i,0)=mult(C(X_i),0)=4p_i$, for all $i\in \Z$.

\hfill

Concerning to the item (b), let $\phi_{ij}: Y_i\to Y_j$ be a bi-Lipschitz homeomorphism
such that $\phi_{ij}(0)=0$.  Let $\widetilde{\phi_{ij}}: \C^{n_i+1}\to \C^{n_j+1}$ 
(resp. $\widetilde{\psi_{ij}}: \C^{n_j+1}\to \C^{n_i+1}$) be a Lipschitz extension of 
$\phi_{ij}$ (resp. $\psi_{ij}=\phi_{ij}^{-1}$) (see \cite{Kirszbraun:1934}, \cite{Mcshane:1934} and \cite{Whitney:1934}).
Let us define $\varphi, \psi: \C^{n_i+1}\times \C^{n_j+1}\to \C^{n_i+1}\times \C^{n_j+1}$ as follows:
$$
\varphi(x,y)=(x-\widetilde{\psi_{ij}}(y+\widetilde{\phi_{ij}}(x)),y+\widetilde{\phi_{ij}}(x))
$$
and
$$
\psi(z,w)=(z+\widetilde{\psi_{ij}}(w), w-\widetilde{\phi_{ij}}(z+\widetilde{\psi_{ij}}(w))).
$$
It easy to verify that $\psi=\varphi^{-1}$ and since $\varphi$ and $\psi$ are composition of Lipschitz maps, they are also Lipschitz maps. Moreover, if $Z_{ij}=Y_i\times \{ 0 \}$ and $\widetilde Z_{ij}=\{ 0\} \times Y_j$, we obtain that $ \varphi(Z_{ij})=\widetilde Z_{ij}$ (see \cite{Sampaio:2016}). Therefore, $(Z_{ij},0)$ and $(\widetilde Z_{ij},0)$ are bi-Lipschitz equivalent, but $mult(Z_{ij},0)=mult(Y_{i},0)$ and $mult(\widetilde Z_{ij},0)=mult(Y_{j},0)$ and, in particular, they have different multiplicity.
\endproof

\hfill


\hfill

{\small

\noindent {\sc Misha Verbitsky\\
            {\sc Instituto Nacional de Matem\'atica Pura e
              Aplicada (IMPA) \\ Estrada Dona Castorina, 110\\
Jardim Bot\^anico, CEP 22460-320\\
Rio de Janeiro, RJ - Brasil}\\
also:\\
{\sc Laboratory of Algebraic Geometry,\\
National Research University Higher School of Economics,\\
Department of Mathematics, 6 Usacheva street, Moscow, Russia.}\\
}}

\noindent {\sc Lev ~Bribrair, Alexandre ~Fernandes and J. Edson ~Sampaio\\ 
			{\sc Departamento de Matem\'atica, Universidade Federal do Cear\'a, \\
Rua Campus do Pici, s/n, Bloco 914\\
Pici, CEP 60440-900\\
Fortaleza-CE, Brazil
}}

\end{document}